\documentclass[11pt,oneside,a4paper]{article}
\usepackage{amssymb}
\usepackage[figuresright]{rotating}
\usepackage{moreverb}
\usepackage{amssymb}
\usepackage{array}
\usepackage{mathrsfs}
\usepackage{amscd}
\usepackage{amsmath}
\usepackage{amsfonts}
\usepackage{amsthm}
\usepackage{psfrag}
\usepackage{textcomp}
\usepackage{bm}
\usepackage{color}
\usepackage{float}
\usepackage[latin1]{inputenc} 	
\usepackage{afterpage}
\usepackage{tabularx}
\usepackage{booktabs}
\usepackage{romannum}
\usepackage{multirow}

\newcommand{\D}[2]{\frac{\partial #1}{\partial #2}}

\newcommand{\ds}	{\displaystyle    }

\newcommand{\p}[1]{\left( #1 \right)}
\newcommand{\cor}[1]{\left[ #1 \right]}
\newcommand{\lla}[1]{\left\{ #1 \right\}}

\newcommand{\om}	{\omega   }

\newcommand{\la}	{\lambda   }

\newcommand{\dt}	{\delta  }
\newcommand{\te}	{\theta   }

\newcommand{\ze}	{\zeta   }

\newcommand{\mbf} {\mathbf  }

\newcommand{\mrm} {\mathrm  }

\newcommand{\bmk}[1] {\bm{\mathfrak{#1}} }

\newcommand{\bmK} {\bm{\mathcal{K}}}

\newcommand{\bmG}{\bm{\mathcal{G}}}

\newcommand{\mD} {\mathcal{D}}

\newcommand{\mG} {\mathcal{G}}

\newcommand{\mbR} {\mathbb{R}}

\newcommand{\mbu}{\mathbf{u}}

\newcommand{\bmC} {\bm{\mathcal{C}}}
\newcommand{\bmu}{\bm{u}}

\newcommand{\bth} {\bm {\theta}}

\newcommand{\ccolumna}[2]{\left\{\begin{array}{c} #1 \\ #2 \end{array}\right\}}

\newcommand{\mmatriz}[4]{\left[\begin{array}{cc} #1 & #2 \\ #3 & #4 \end{array}\right]}

\newcommand{\mbG}{\mathbf{G}}

\usepackage{authblk}
\usepackage[lmargin=2cm,rmargin=2.0cm,tmargin=2cm,bmargin=2cm]{geometry}

\title{Critical damping in nonviscously damped linear systems}

\author[1]{Mario L\'azaro}

\affil[1]{Dep. of Continuum Mechanics and Theory of Structures.
         Universit\`at Politecnica de Val\`encia.
				 46022, Valencia, Spain}

\begin{document}

\maketitle

\begin{abstract}
In structural dynamics, energy dissipative mechanisms with non-viscous damping are characterized by their dependence on the time-history of the response velocity, mathematically represented by convolution integrals involving hereditary functions. Combination of damping parameters in the dissipative model can lead the system to be overdamped in some (or all) modes. In the domain of the damping parameters, the thresholds between induced oscillatory and non--oscillatory motion are called critical damping surfaces (or manifolds, since we can have a lot of parameters). In this paper a general method to obtain critical damping surfaces for nonviscously damped systems is proposed. The approach is based on transforming the algebraic equations which defined implicitly the critical curves into a system of differential equations. The derivations are validated with three numerical methods covering single and multiple degree of freedom systems. 
\end{abstract}

%
%
%

\section{Introduction}


In this paper, nonviscously damped linear systems are under consideration. Nonviscous (also named by viscoelastic) materials have been widely used for vibrating control in mechanical, aerospace, automotive and civil engineering applications. This paper deals precisely with those applications where vibrations are tried to be disappeared, that is,  designing of damping devices which are able to avoid oscillatory motion at dynamical systems. In nonviscous models, damping forces are assumed to be dependent on the history of the response velocity via  kernel time functions. As far as the motion equations concerned, this fact is represented by convolution integrals involving the velocities of the degrees--of--freedom (dof) and affected by the hereditary kernels. Denoting by $\bmk{u}(t) \in \mathbb{R}^n$ to the array with the degrees of freedom of the system, this vector verifies the dynamic equilibrium equations which in turn has an integro-differential form
\begin{equation}
		\mbf{M}\ddot{\bmk{u}} + \int_{-\infty}^t \bm{\mathcal{G}}(t-\tau) \, \dot{\bmk{u}}\ \mrm{d}\tau + \mbf{K}\bmk{u} = \bmk{f}(t)
		\label{eq001}
\end{equation}		
where $\mbf{M}, \, \mbf{K}  \in \mathbb{R}^{n\times n}$ are the mass and stiffness matrices assumed to be positive definite and positive semidefinite, respectively; $\bm{\mathcal{G}}(t) \in \mathbb{R}^{n\times n}$ is the nonviscous damping matrix in the time domain, assumed symmetric, which satisfies the necessary conditions of Golla and Hughes~\cite{Golla-1985} for a strictly dissipative behavior. As known, the viscous damping is just a particular case of Eq.~\eqref{eq001} with $\bm{\mathcal{G}}(t) \equiv \mbf{C} \, \delta(t)$, where $\mbf{C}$ is the viscous damping matrix and $\delta(t)$ the Dirac's delta function. The time--domain system of motion equations are then reduced to the well known expressions
\begin{equation}
		\mbf{M}\ddot{\bmk{u}} + \mbf{C} \dot{\bmk{u}} + \mbf{K}\bmk{u} =\bmk{f}(t)
		\label{eq001b}
\end{equation}		
Considering now the free motion case $\bmk{f}(t) \equiv \mbf{0}$ in Eq.~\eqref{eq001}, we test exponential solutions of the type $\bmk{u}(t) =\bmu \, e^{st}$, with $\bmu$ and $s$ to be found. Then, the classical nonlinear eigenvalue problem associated to viscoelastic vibrating structures yields
\begin{equation}
		\left[s^2\mbf{M} + s \mbf{G}(s) + \mbf{K} \right]\bmu \equiv \mbf{D}(s) \, \bmu = \mbf{0}
		\label{eq002}
\end{equation}		
where $\mbf{G}(s) = \mathcal{L}\{ \bm{\mathcal{G}}(t) \} \in \mathbb{C}^{N\times N}$ is the damping matrix in the Laplace domain and $\mbf{D}(s)$ is the dynamical stiffness matrix or transcendental matrix. \\

Response of Eqs.~\eqref{eq001} is closely related to the eigensolutions of the eigenvalue problem~\eqref{eq002}. Adhikari~\cite{Adhikari-2002c} derived modal relationships and closed form expressions for the transfer function in the Laplace domain. Due to the non--linearity, induced by a frequency-dependent damping matrix, the search eigensolutions is in general much more computationally expensive than that of classical viscous damping~\cite{Adhikari-2003c}. A survey of the different viscoelastic models can be found on the references~\cite{Adhikari-2010b,Lazaro-2012} although as far as this paper concerned we work with hereditary damping models based on exponential kernels~\cite{Biot-1955}.\\

Since the present paper is devoted on critical damping and this field has been deeply studied in  the bibliography for viscously damped systems, we consider relevant to review the main works on it. Duffin~\cite{Duffin-1955} defined  an overdamped system in terms of the quadratic forms of the coefficient matrices.  Nicholson~\cite{Nicholson-1978} obtained eigenvalue bounds for free vibration of damped linear systems. Based on these bounds, sufficient condition for subcritical damping were derived. M\"uller~\cite{Muller-1979} characterized an underdamped system in similar terms to Duffin's work deriving a sufficient condition expressed as function of the definiteness of the system matrices. Inman and Andry~\cite{Inman-1980} proposed sufficient conditions  for underdamped, overdamped and critically damped motions in terms of the definiteness of the system matrices. These conditions are valid for classically damped systems although Inman and Andry shown that they also could work for non--classical systems. Inman and Orabi~\cite{Inman-1983} and Gray and Andry~\cite{Gray-1982} proposed more efficient method for computing the critical damping condition. However, Barkwell and Lancaster~\cite{Barkwell-1992} pointed out some defficiencies in the Inman and Andry criterion of ref.~\cite{Inman-1980} presenting a counterexample and they provided some reasons explaining why this criterion had been usually adopted to check criticality in damped systems. Additionally, Barkwell and Lancaster~\cite{Barkwell-1992} obtained necessary and sufficient conditions for overdamping in gyroscopic vibrating sytems. Bhaskar~\cite{Bhaskar-1997} presented a more complete overdamping condition which somehow corrected that of Inman and Andry~\cite{Inman-1980} giving a generalization to a class of non--conservative systems. Beskos and Boley~\cite{Beskos-1980} established conditions for finding critical damping surfaces from the determinant of the system and its derivative. They proved that a critically damped eigenvalue was simultaneously root of the characteristic equation and its derivative, something that can be used to detect critical damping surfaces. In the work~\cite{Beskos-1981} the same authors studied conditions for critical damping in continuous systems. Beskou and Beskos~\cite{Beskos-2002} presented an approximate method computationally efficient to find critical damping surfaces separating overdamping or partially overdamping regions from those of underdamping for viscously damped systems. 	\\

As far as nonviscous systems concerned, research on critical damping has not been as exhaustive as that of viscous damping. Mainly, investigations have been focused on single degree--of--freedom systems on the discussion of the type of response attending at the damping parameters of a single--exponential hereditary kernel.  Muravyov and Hutton \cite{Muravyov-1998b} and Adhikari~\cite{Adhikari-2005} analyzed the conditions under which single degree--of--freedom nonviscously damped systems by one exponential kernel becomes critically damped. He carried out an exhaustive analysis of the roots nature of the resulting third order characteristic polynomial. Adhikari~\cite{Adhikari-2008} studied the dynamic response of nonviscously damped oscillators and discussed the effect of the damping parameters on the frequency response function. M\"uller~\cite{Muller-2005} performed a detailed analysis on the nature of the eigenmotions of a single degree of freedom Zener 3--parameter viscoelastic model. Muravyov~\cite{Muravyov-1998b} obtained closed--form solutions for forced nonviscoulsy damped beams studying the conditions for overdamping or underdamping time response. As known, viscoelastic systems modeled by hereditary exponential functions are characterized by having extra real overdamped modes associated to those kernels. As far as these type of modes concerned the references~\cite{Lazaro-2014a,Lazaro-2015d,Lazaro-2016c} provide a mathematical characterization and some numerical methods to their evaluation.L\'azaro~\cite{Lazaro-2012} observed that certain recursive method to obtain eigenvalues in proportionally damped viscoelastic systems always converges under a linear rate except just in the critical surfaces where the scheme is underlineal. \\


In this paper, critical damping surfaces of nonviscously damped linear systems are presented. Critical damping is refereed to the set of damping parameters within the threshold between induced oscillatory and non--oscillatory motion (for all or for some modes). The general procedure to extract these manifolds in the domain of the damping parameters is to eliminate a parameter of a system of two algebraical equations. Encouraged by the fact that this elimination is not possible for polynomials with order greater than four, a new method to construct critical curves is developed. This method is based on to transform the algebraical equations into two ordinary differential equations. The method is validated with three numerical examples. The two first are devoted on single degree--of--freedom systems with one and two hereditary exponential kernesl,  respectively. The application of the current approach for multiple degree--of--freedom sytems is presented in the third example.

\section{Conditions of criticality in terms of determinant of the system}
\label{ProposedMethod}

In this section we will extend the main results derived by Beskos and Boley~\cite{Beskos-1980}  on critical viscous damping to nonviscously damped systems. In order to establish the basis of our work, we will describe the type of damping model adopted in its most general form. We will consider a damping matrix based on hereditary Biot's exponential kernels. Mathematically, this model adopts the following form in time  and in frequency domain 
\begin{equation}
\bmG(t) = \sum_{k=1}^N  \mbf{C}_k \, \mu_k \,  e^{-\mu_k t}  \ ,  \quad
\mbG(s) = \mathcal{L}\{\bmG(t)\} = \sum_{k=1}^N \frac{\mu_k}{s + \mu_k}  \mbf{C}_k
\label{eq004b}
\end{equation}
where $\mu_k>0$, $1\leq k \leq N$ represent the relaxation or also called nonviscous coefficients and $\mbf{C}_k \in \mathbb{R}^{n \times n}$ are the (symmetric) matrices of the limit viscous damping model, defined as the limit
\begin{equation}
		 \sum_{k=1}^N  \mbf{C}_k = \lim_{\mu_1 \ldots \mu_N \to \infty} \mbG(s) 
		\label{eq005}
\end{equation}
Coefficients $\mu_k$ control the time and frequency dependence of the damping model while the spatial location and the level of damping are controlled by coefficients within matrices $\mbf{C}_k$. It is straightforward that the following relationships hold
\begin{equation}
		\sum_{k=1}^N  \mbf{C}_k = \int_{0}^\infty \bmG(t) \, dt  = \mbG(0)
		\label{eq006}
\end{equation} 
Henceforth, we will consider the damping matrix, and for extension the transcendental matrix, as depending not only on the frequency via $s$, but also  on  set of parameters controlling the dissipative behavior. In the most general case, let us see that the symmetric damping model presented in \eqref{eq004b} depends on $p_{\text{max}} = N + Nn(n+1)/2$ parameters. Indeed, $N$ nonviscous coefficients $\mu_1,\ldots,\mu_N$ plus $n(n+1)/2$ possible independent entrees in every symmetric matrix $\mbf{C}_k$, with $1\leq k \leq N$. Thus, the complete set of parameters can be listed as
\begin{equation}
\mu_1,\ldots,\mu_N,C_{111},\ldots,C_{1nn},\ldots,C_{N11},\ldots,C_{Nnn}
\label{eq009}
\end{equation}
where $C_{kij}=C_{kji}$ is the entree $ij$ of matrix $\mbf{C}_k$. Real applications depend in general on much less parameters, say $p << p_{\text{max}}$. In the sake of clarity, we will denote by $\bth = \lla{\te_1,\ldots,\te_p}$ the set of independent damping parameters and consequently we can write the damping matrix as $\mbf{G}(s,\bth)$.\\

According to the said above, we can denote as $\mD(s,\bth) = \det \cor{\mbf{D}(s)}$ the determinant associated to the nonlinear eigenvalue problem~\eqref{eq002}. Eigenvalues are then roots of the equation
\begin{equation}
		\mD(s,\bth) =  \det \cor{ s^2\mbf{M} + s \sum_{k=1}^N \frac{\mu_k}{s + \mu_k}  \mbf{C}_k + \mbf{K} } = 0
		\label{eq007}
\end{equation}
Attending to the values of $\bth \in \mbR^p$, the algebraic structure of the spectrum of this problem can vary. If the level of damping induced by matrix $\mbf{G}(s,\bth)$ is light, we will have $2n$ complex eigenvalues with oscillatory nature and $r$ real eigenvalues  with non--oscillatory nature and associated to the nonviscous hereditary kernels (hence  they are also usually named as nonviscous eigenvalues). Furthermore, the total number of these real eigenvalues is~\cite{Adhikari-2003c} 
	\begin{equation}
	r = r_1 +  \cdots + r_N = \sum_{j=1}^N \text{rank} (\mbf{C}_j)
	\label{eq008}
	\end{equation}
As long as there exist $2n$ complex eigenvalues and $r$ nonviscous eigenvalues, we will say that the  system is completely underdamped. As the damping level increases, the real part of eigenvalues (not necessary all) becomes higher (in absolute value) and the imaginary part decreases. For certain value of the damping parameters a conjugate--complex pair could merge into a double real negative root. The set of damping parameters is said then to be on a critical surface, which in turn represents the threshold between underdamping and overdamping. If oscillatory modes coexist with those non--oscillatory, then we say the system is partially overdamped (or mixed overdamping). The system is said to be completely overdamped if all modes are so. For mixed or complete overdamping, some (or maybe all) of the roots of~\eqref{eq007} are negative real numbers, say $s = \la$, with $\la < 0$ so that
\begin{equation}
\mD(\la,\bth) =  \det \cor{ \la^2\mbf{M} + \la \mbf{G}(\la,\bth)  + \mbf{K} } = 0
\label{eq007b}
\end{equation}
For each value of $\la$, Eq.~\eqref{eq007b} defines a $p$--dimensional surface in the space where the parameters array $\bth$ can take values. Since for light damping we have initially $n$ pairs of conjugate--complex eigenvalues, we will have as much as $n$ critical surfaces because, as Beskos and Boley~\cite{Beskos-1980} point out: ``there are at most as many partial critical damping possibilities as the number of the pairs in~\eqref{eq007} of roots $s$ with zero imaginary part''. The mathematical principle which characterizes a critical damping surfaces can be extrapolated to non--vicous damping and therefore these ones can be found imposing a minimum among all possible values of $\la$ in Eq.~\eqref{eq007b}, that is
\begin{equation}
\D{}{\la}\mD(\la,\bth) = \D{}{\la}\det \cor{\mbf{D}(\la,\bth)}=0
\label{eq011}
\end{equation}
This condition is consistent with the fact that a critical root is double just under critical condition. Therefore, Eqs.~\eqref{eq007} and~\eqref{eq011} define a set of critical surfaces resulting after eliminating parameter $\la$ from both equations. This process, although well defined from a theoretical point of view, can only be carried out if an analytical closed--form of $\mD(\la,\bth)$ is provided something that only occurs for small to moderately sized systems. For nonvisocusly damped systems, this procedure has not been used yet, to the authors knowledge. Instead, they have been found for single degree--of--freedom systems and for $N=1$ kernels since this particular problem leads to a three order polynomial, which as known allows radicals--based analytical solution (Cardano's formulas). Additionally, in the present paper we also attempt to improve the numerical evaluation of critical damping surfaces proposing a numerical method which will be described in the following paragraphs.\\

Numerical evaluation of critical surfaces consists in solving Eqs.~\eqref{eq007} and~\eqref{eq011} simultaneously for a prefixed range of values of damping parameters. This process  becomes in general computationally inefficient since for each value of the prefixed parameters, a system of two non--linear  equations must be solved. Attempting to improve the numerical procedure for constructing critical surfaces we propose a method valid to build critical curves formed by two parameters, assuming as fixed the rest of them. The method is able to find critical overdamped regions in two dimensional cross sections of the $p$--dimensional real critical manifolds. Thus, from the complete set of parameters $\bth = \{\te_1,\ldots,\te_p\}$, we chose two of them, which will be named as design parameters. Without loss of generality, we can take $\te_1$ and $\te_2$ while the rest of parameters  remain fixed, say $\te_{30},\ldots,\te_{p0}$. The  challenge is to draw the critical damping curves in the plane $(\te_1,\te_2)$. For a sake of clarity in the notation, we will denote by $p = \te_1$ and $q = \te_2$ and will assume then that the critical curve(s) are functions of the form $q = q(p)$. For each value of $p$, both equations
\begin{equation}
\mD(\la,p,q) = 0 \quad , \qquad \D{}{\la} \mD(\la,p,q) = 0
\label{eq012}
\end{equation}
allow to find a pair ($\la, q$) (or several, since $\la$ is within a polynomial). Let us consider a point $p_0$ for which $q_0$ and $\la_0$ are solutions of Eqs.~\eqref{eq012} and let us assume around $p=p_0$ the functions $q(p)$ and $\la(p)$ exist. The three numbers $(p_0,q_0,\la_0)$ form a initial point of the proposed approach. The derivatives 
$\la'(p)  =\mathrm{d}\la / \mathrm{d}p$ and $q'(p) = =\mathrm{d}q / \mathrm{d}p$ 
can be evaluated just  applying the chain rule in Eqs.~\eqref{eq012}. Indeed, 
\begin{eqnarray}
\mD_{,\la} \, \la'(p) + \mD_{,q} \, q'(p) + \mD_{,p} &=& 0 \\
\mD_{,\la \la} \, \la'(p) + \mD_{,\la q} \, q'(p) + \mD_{,\la p} &=& 0
\label{eq013}
\end{eqnarray}
where subscripts denoting partial derivatives. Since from Eq.~\eqref{eq012}, we have $\mD_{,\la} = 0$,   then we can solve for $\la'(p)$ and $q'(p)$
\begin{eqnarray}
q'(p) &=& - \frac{\mD_{,p}}{\mD_{,q}} \nonumber \\
\la'(p) &=& \frac{ \mD_{, p} \mD_{,\la q} }{ \mD_{,q} \mD_{,\la \la}} - \frac{ \mD_{,\la p} }{\mD_{,\la \la}}
\label{eq014}
\end{eqnarray}
These two equations form a system of two ordinary differential equations whose solutions are be well defined provided that the derivatives $\mD_{,\la \la}$ and $\mD_{,q}$ do not vanish at $(p_0,q_0,\la_0)$. Existence of the critical curve beyond a close interval around the initial point will be subordinate to the existence of those derivatives along the curve. Critical curves arise now  as the numerical solution of a system of ordinary differential equations, for which Runge--Kutta based methods can be used. Before, the method requires solving the system of two algebraical equations~\eqref{eq012} and two unkowns, say $\la_0, \, q_0$, which in general results in several solutions because of the polynomial form. Pairs $(\la_0,q_0)$ both reals and verifying $\la_0 < 0$ and  $q_0 \geq 0$ (we will assume a positive range for parameters) will be appropriate solutions lying on a critical surface. Taking derivatives repeatedly respect $p$ in~\eqref{eq013} also lead us to obtain higher order derivatives, allowing to find the Taylor expansion of the critical curve around $p=p_0$. This procedure is used to find an approximation of a critical curve in the one--kernel single--dof numerical example.

\section{Numerical examples}

\subsection{Single degree of freedom systems, $N=1$ exponential kernel}

%
\begin{figure}[ht]
\begin{center}
\input{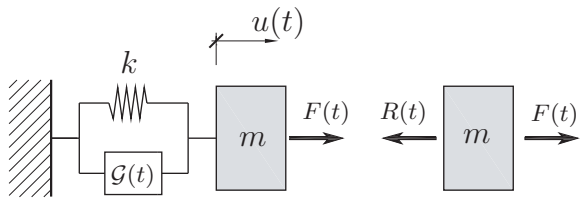}
\caption{A single degree--of--freedom viscoelastic oscillator}
\label{fig05}
\end{center}
\end{figure}

First, we will consider the single dof nonviscous system with one hereditary exponential kernel. The dof represents the displacement of certain mas $m$ attached to ground by the  viscoelastic constraint. Fig.~\ref{fig05} shows the schematic configuration {\em mass--spring--viscoelastic damper} and the corresponding free body diagram. Hence, the internal force  is related to the displacement by
\begin{equation}
		R(t) = \int_{-\infty}^t \mathcal{G}(t-\tau) \dot{u}(\tau) \mrm{d}\tau + k u(t) 
		\label{eq015}
\end{equation}		
$k$ is the constant of the linear--elastic spring and $\mathcal{G}(t)$ is the dissipative kernel or damping function with the general form, both in time and frequency domain
\begin{equation}
\mathcal{G}(t) = c \, \mu \, e^{-\mu t} \ , \quad G(s) = \frac{\mu \, c}{s + \mu}
\label{eq016}
\end{equation}
where $\mu$ and $c$ are respectively is the nonviscous and the viscous coefficients. The free motion equation can be deduced from the dynamic equilibrium for $F(t) \equiv 0$
\begin{equation}
			m\ddot{u} + \int_{-\infty}^t \mathcal{G}(t-\tau) \dot{u}(\tau) \mrm{d}\tau + k u(t) = 0
			\label{eq017}
\end{equation}		
And the associated characteristic equation  
\begin{equation}
			m s^2 + s G(s) + k =  m s^2 + s \frac{\mu \, c}{s + \mu} + k = 0
			\label{eq018}
\end{equation}		
It is quite appropriate to board this problem using dimensionless variables in order to compare with existing results presented in the bibliography. Thus, we define the following non--dimensional variables
\begin{equation}
x = \frac{s}{\om_n} \ , \quad \nu = \frac{\om_n}{\mu} \ , \quad \ze = \frac{c}{2 m \om_n}
\label{eq019}
\end{equation}
where $\om_n= \sqrt{k/m}$ is he natural frequency of the undamped system. After straight operations and multiplying Eq.~\eqref{eq018} by the denominator of the damping function we obtain the characteristic equation in non--dimensional form as the third order polynomial
\begin{equation}
\mD(x,\nu,\ze) = (1 + \nu x) ( x^2 + 1) + 2 x \ze = \nu x^3 + x^2 + (\nu+2\ze)x + 1 = 0
\label{eq020}
\end{equation}
As known, the three roots are available as function of the coefficients so that a detailed discussion of the nature of the three roots can be addressed as function of the values of $\ze>0$ and $\nu>0$. This work was carried out by Adhikari in the references~\cite{Adhikari-2005,Adhikari-2008} where closed form expressions of the critical curves enclosing the overdamped region were derived. For the  sake of our exposition we consider interesting transcript here the Adhikari's results of the critical curves, since later we will present also approximations of them. Thus, the overdamped region can be defined as the set 
\begin{equation}
\lla{(\ze,\nu) \in \mbR^+ : \ze_L(\nu) \leq \ze \leq \ze_U(\nu)} 
\label{eq021}
\end{equation}
where the critical damping curves $\ze_L(\nu),\,\ze_U(\nu)$ are
\begin{eqnarray}
		\ze_L(\nu) &=& \frac{1}{24 \nu} \cor{1 - 12\nu^2 + 2 \sqrt{1 + 216\nu^2} + \cos \p{\frac{4\pi + \te}{3}}} \nonumber \\
		\ze_U(\nu) &=& \frac{1}{24 \nu} \cor{1 - 12\nu^2 + 2 \sqrt{1 + 216\nu^2} + \cos \p{\frac{\te}{3}}} \label{eq022}
\end{eqnarray}
with
\begin{equation}
\te = \arccos \cor{ -\frac{5832\nu^4 + 540 \nu^2 - 1}{\p{1+216\nu^2}^{3/2}}}
\label{eq023}
\end{equation}

Let us apply the proposed method to find critical damping curves based on the solution of the system of differential equations~\eqref{eq014}. According to the theoretical derivations, the critical surfaces arises from eliminating $x$ from the two following equations
\begin{eqnarray}
\mD(x,\nu,\ze) &=&  \nu x^3 + x^2 + (\nu+2\ze)x + 1 = 0 \nonumber \\
\D{\mD}{x}     &=& 3\nu x^2 + 2x  + \nu + 2 \ze     = 0
\label{eq024}
\end{eqnarray}
From the second equation we can obtain the two roots $x_{1,2} = (-1 \pm \sqrt{-6 \zeta  \nu -3 \nu ^2+1})/3\nu$ and then plug them into the first one. After some simplifications, we obtain the critical surfaces in implicit form
\begin{equation}
8 \zeta ^3 \nu +12 \zeta ^2 \nu ^2-\zeta ^2+6 \zeta  \nu ^3-10 \zeta  \nu +\nu ^4+2 \nu ^2+1=0
\label{eq025}
\end{equation}
which coincides with the third order polynomial obtained by Adhikari~\cite{Adhikari-2005}. \\

\begin{figure*}[htb]%
\begin{center}
	\input{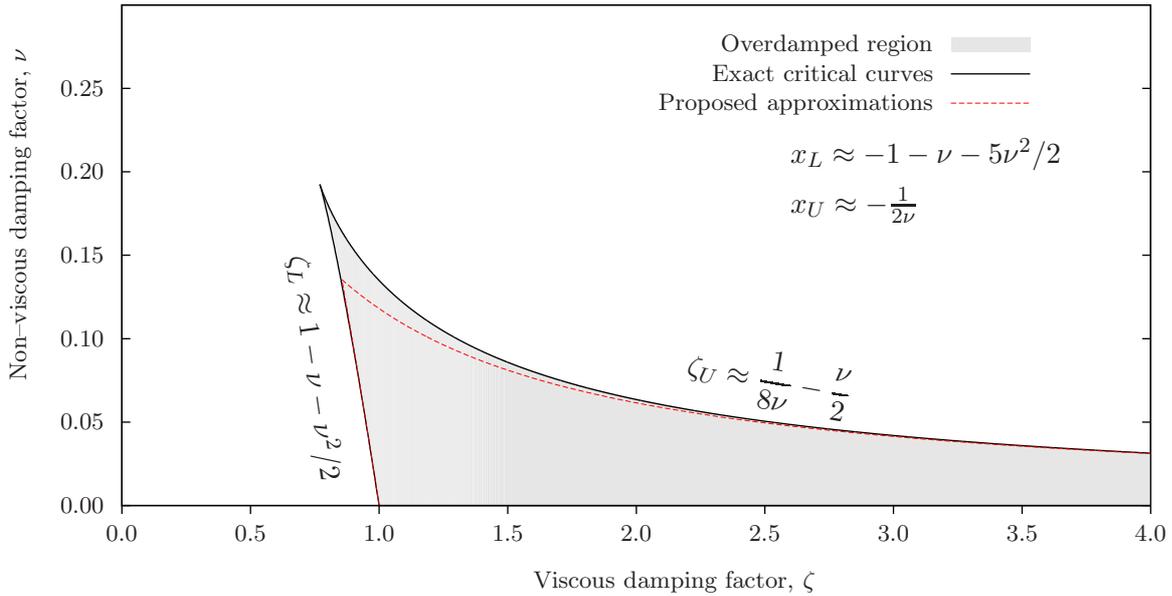}  \\	
\caption{Example 1: single degree--of--freedom with $N=1$ exponential kerkel. Exact and approximate (proposed) overdamped region. $x_L$ and $x_U$ are the overdamped critical eigenvalues along the critical (approximate) curves}%
\label{fig01}%
\end{center}
\end{figure*}

We will attempt now to obtain curves of the form $\ze = \ze(\nu)$, therefore our independent parameter is $p=\nu$ and the dependent variable is $q = \ze$. We need to find the partial derivatives of $\mD(x,\nu,\ze)$ and $\mD_{,x}(x,\nu,\ze)$ respect to $x$, $\ze$ and $\nu$, obtaining
\begin{eqnarray}
\mD_{,\nu} &=& x + x^3 \quad , \qquad  \mD_{,\ze} = 2x \nonumber \\
\mD_{,x\nu} &=& 1 + 3x^2 \quad , \qquad  \mD_{,x\ze} = 2 \quad , \qquad  \mD_{,xx} = 2 + 6 \nu x
\label{eq026}
\end{eqnarray}
After some math, the two differential equations are set as
\begin{eqnarray}
\ze'(\nu) &=& - \frac{2}{1 + x^2} \nonumber \\
x'(\nu) &=&  \frac{2x^2}{(1+x^2)(1+3\nu x)}
\label{eq027}
\end{eqnarray}
These equations must be completed with initial conditions. Taking $\ze_0 = 1$ equations~\eqref{eq024} can be solved obtaining four pairs of roots
\begin{equation}
(x = -1, \ \nu = 0)  \quad ; \quad   (x = -3.38298, \ \nu =0.134884)  \quad ; \quad  (x = 0.191 \pm 0.508 i, \ \nu = -3.067 \pm 2.327i) 
\label{eq028}
\end{equation}
Only real solutions with $x<0, \, \nu\geq 0$ are of interest as initial values. The first pair results in the initial values $\nu_0=0,\ze_0=1,x_0=-1$ of the critical curve $\ze_L(\nu)$ while the second pair $\nu_0=0.134884,\ze_0=1,x_0=-3.38298$ gives as a result $\ze_U(\nu)$. Both curves have been plotted in Fig.~\eqref{fig01}. Results fit perfectly with those of exact solutions of Eqs.~\eqref{eq022}.
Although from an analytical point of view, the problem of determining critical curves in this case is solved, we will go further proposing two closed--form approximations of $\ze_L(\nu)$ and  $\ze_U(\nu)$. We consider these derivations of interest on one hand,  by their simplicity respect to those of the exact expressions. And, on  the other hand,  by the procedure to deduce them. \\

Encouraged by the fact that the initial point of the critical curve $\ze_L(\nu)$ is as simple as $\ze_L(0)=1$ and also by its regularity and low curvature (information already known since the exact result is available in Fig.~\ref{fig01}, we think that the Taylor series expansion around $\nu_0=0$ can provide accurate results and in turn simple in form. Indeed, the first derivative can be determined just from Eq.~\eqref{eq027} for $x_0=-1$ and $\nu_0=0$, resulting 
$$\ze'_L(0) = -1 \ , \quad x'(0)=-1$$
Now, taking again derivatives respect to $\nu$ in Eqs.~\eqref{eq026} and after some operations, second derivatives $\ze_L''(0)$ and $x''(0)$ can be found, so that  
$$\ze''_L(0) = -1 \ , \quad x''(0)=-5$$
Hence, Taylor series expansions up to the second order of the critical damping curve $\ze_L(\nu)$, and its associated critical eigenvalue $x_L(\nu)$ are then

\begin{eqnarray}
\ze_L(\nu) &\approx & 1 - \nu - \nu^2/2 \label{eq029a} \\
x_L(\nu) & \approx & -1 -\nu -5 \nu^2/2 \label{eq029b}
\end{eqnarray}

Similar procedure could be followed to find a Taylor based approximation around upper critical curve $\ze_U$, however, this function presents higher changes of curvatures and a wider domain of $\ze$ (in fact an infinite range). In is expected that a polynomial based approximation only will work around the initial point and of course it will not be able to represent the asymptotic behavior. To undertake this approach a recent result on asymptotic behavior of polynomial roots proposed by L\'azaro et al.~\cite{Lazaro-2015a} will be used. Given a polynomial 
\begin{equation}
a_0 + a_1 X + \cdots +    a_{n-2} X^{n-2} + a_{n-1} X^{n-1} + X^{n} 
\label{eq030}
\end{equation}
then the numbers (called polynomial pivots)
\begin{equation}
- \frac{a_{n-1}}{2} \pm \sqrt{ \p{\frac{a_{n-1}}{2}}^2 - a_{n-2}}
\label{eq031}
\end{equation}
present the property of lying close to one (or two) roots provided that they are not relatively smaller than the rest of the polynomial coefficients. The exact mathematical conditions describing this statement are given in form of several theorems in the reference~\cite{Lazaro-2015a}. We check if the so defined pivots of the third order polynomial $\mD(x,\nu,\ze)$ can give us valuable information respect to the nature of the roots. Thus, the  pivots are
\begin{equation}
- \frac{a_{2}}{2} \pm \sqrt{ \p{\frac{a_{2}}{2}}^2 - a_{1}} = - \frac{1}{2\nu} \pm \sqrt{\frac{1}{4\nu^2}  - \frac{2\ze}{\nu} - 1 }
\label{eq032}
\end{equation}
Looking for critical damping curves, we know that along them, the roots have double multiplicity (double roots). Therefore, if we admint that the pivots are close to two roots of the problem and we force the discriminant of Eq.~\eqref{eq032} to be zero, then the obtained root will be double and we will lie on a critical damping curve. Hence, to vanish the discriminant results in the approximation of the upper curve $\ze_U(\nu)$. Indeed,
\begin{equation}
\frac{1}{4\nu^2}  - \frac{2\ze}{\nu} - 1 = 0 \quad \to \quad \ze_U(\nu) \approx \frac{1}{8 \nu } -\frac{\nu }{2}
\label{eq033}
\end{equation}
The associated damped eigenvalue can be approximated by
\begin{equation}
x_U(\nu) \approx - \frac{1}{2\nu} 
\label{eq034}
\end{equation}
According to the results of~\cite{Lazaro-2015a} the bigger the pivots the closer to the roots of the polynomial. Therefore, we can predict that the lower the nonviscous parameter $\nu$ the more accurate the results, as can be appreciated in the Fig.~\ref{fig01}.

\subsection{Single degree of freedom systems, $N=2$ exponential kernels}

Now we will attempt to find the critical damping surfaces of a single dof system with $N=2$ hereditary kernels. The approach can easily be extrapolated to the general case of $N$ kernels. According to Eq.~\eqref{eq004b}, the damping function is
\begin{equation}
\mG(t) = c_1 \, \mu_1 \,  e^{-\mu_1 t} + c_2 \, \mu_2 \,  e^{-\mu_2 t}  \ ,  \quad
G(s) = \mathcal{L}\{\mG(t)\} = \frac{\mu_1 \, c_1}{s + \mu_1}  + \frac{\mu_2 \, c_2}{s + \mu_2} 
\label{eq035}
\end{equation}
And the characteristic equation yields
\begin{equation}
			m s^2 + s G(s) + k =  m s^2 + s  \p{\frac{\mu_1 \, c_1}{s + \mu_1}  + \frac{\mu_2 \, c_2}{s + \mu_2}}  + k = 0
			\label{eq036}
\end{equation}		
As noticed, the dissipative model has four parameters, $c_1,c_2,\mu_1,\mu_2$. Our procedure allows to draw critical curves of two parameters, hence before finding the solution of the proposed differential equations two parameters must be fixed. For a sake of the solution representation, we fusion the damping coefficients into only one, considering the case $c_1=c_2=c$. Let us define the following dimensionless parameters
\begin{equation}
x = \frac{s}{\om_n} \ , \quad 
\nu_1 = \frac{\om_n}{\mu_1} \ , \quad 
\nu_2 = \frac{\om_n}{\mu_2} \ , \quad 
\ze = \frac{c}{m \om_n}
\label{eq037}
\end{equation}
Where $\om_n = \sqrt{k/m}$ is the natural frequency of the  system. The Eq.~\eqref{eq036} can be expressed now in dimensionless form
\begin{equation}
			x^2 + x \, \ze   \p{\frac{1}{1 + \nu_1 \, x}  + \frac{1}{1 + \nu_2 \, x}}  + 1 = 0
			\label{eq039}
\end{equation}		
Note that $\nu_1=\nu_2=0$ yields the particular case of viscous damping leading to $\ze_{cr}=1$. Multiplying Eq.~\eqref{eq039} by $(1 + \nu_1 x)(1 + \nu_2 x)$ we transform the characteristic equation into a four order polynomial equation
\begin{equation}
\mD(x,\nu_1,\nu_2,\ze) = (1 + \nu_1 x)(1 + \nu_2 x) ( x^2 + 1) + x \ze (2 + \nu_1 x + \nu_2 x)  = 0
\label{eq038}
\end{equation}
which together with 
\begin{equation}
\mD_{,x}(x,\nu_1,\nu_2,\ze) = 2 \ze + 2x + 2x \, \nu_1 \nu_2 (1 + 2x^3) + (1 + 2x\ze + 3 x^3)(\nu_1 + \nu_2)
\label{eq040}
\end{equation}
allow to find the critical surfaces after eliminating $x$. Observe that the  last equation is a third order polynomial, therefore the Cardano formulas can be used to obtain the three roots. Plugging them into Eq.~\eqref{eq038} would lead the exact critical curves. These derivations will not be carried out here since the resulting expressions would be hardly handled and together with the difficulty to follow properly the exposition. On the other hand, they can easily be programmed in a symbolic software and results be compared to those of the present method. \\

According to the proposed approach, the critical surfaces are defined in terms of two parameters, leaving fixed the rest. For the  current example we will consider critical curves in the plane $(\ze,\nu_2)$ (i.e. functions $\nu_2 = f(\ze)$ with $\nu_1$ fixed) and also in the plane $(\nu_1,\nu_2)$ (i.e. functions $\nu_2 = f(\nu_1)$ with $\ze$ fixed). Two systems of differential equations must be assembled, one involving the functions $x(\ze), \nu_2(\ze)$ and the other one the functions $x(\nu_1), \nu_2(\nu_1)$. From Eq.~\eqref{eq013} the two problems can be written in matrix form as
\begin{itemize}
	\item Critical curves in plane $(\ze,\nu_2)$. Parameter $\nu_1$ fixed. 
\begin{equation}
\mmatriz{0}{\mD_{,\nu_2}}{\mD_{,xx}}{\mD_{,x\nu_2}} 
\ccolumna{x'(\ze)}{\nu_2'(\ze)} = 
- \ccolumna{\mD_{,\ze}}{\mD_{,x \ze}}
\ , \quad
\ccolumna{x(\ze_{0})}{\nu_2(\ze_{0})} = \ccolumna{x_0}{\nu_{20} } 
\label{eq041}
\end{equation}
	\item Critical curves in plane $(\nu_1,\nu_2)$. Parameter $\ze$ fixed. 
\begin{equation}
\mmatriz{0}{\mD_{,\nu_2}}{\mD_{,xx}}{\mD_{,x\nu_2}} 
\ccolumna{x'(\nu_1)}{\nu_2'(\nu_1)} = 
- \ccolumna{\mD_{,\nu_1}}{\mD_{,x \nu_1}}
\ , \quad
\ccolumna{x(\nu_{10})}{\nu_2(\nu_{10})} = \ccolumna{x_0}{\nu_{20} } 
\label{eq042}
\end{equation}
\end{itemize}
where
\begin{align*}
\mD_{,\nu_1}&= x \cor{1 + x (x + \ze + \nu_2 + x^2 \nu_2)}     &  
\mD_{,\nu_2} &= x \cor{1 + x (x + \ze + \nu_1 + x^2 \nu_1)}   &
\mD_{,\ze}&= x \cor{2 + x (\nu_1 + \nu_2)} \\
\mD_{,x\nu_1}&=1 + x \cor{2 (\ze + \nu_2) + x (3 + 4 x \nu_2)}      &
\mD_{,x\nu_2} &= 1 + x \cor{2 (\ze + \nu_1) + x (3 + 4 x \nu_1)}              &
\mD_{,x\ze}&=2 \cor{1 + x (\nu_1 + \nu_2)}
\end{align*}
The initial conditions come from solving Eqs.~\eqref{eq038} and~\eqref{eq040} for prescribed values of two of the parameters. Table~\ref{tab01} lists a complete set of initial values which allow to address the solution of the differential equations. From the solution of the aforementioned algebraic equations one can find simultaneously two different initial conditions. Thus, for instance, for the case $\nu_{10}=0.00$ and $\ze_0=1.00$, among other complex solutions we find $(\nu_{20},x_0)=(0,-1) $ and $ (\nu_{20},x_0)=(0.1916,-2.7693)$ as the initial conditions of the two first curves shown in Table~\ref{tab01}. Notice that this fact also takes place for other cases in the table, always if we somewhat are fortunate in our election of the prescribed pair of parameters. Otherwise, we could not obtain satisfactory solutions, for example taking $\nu_{10}=0.25$ and $\ze_0 = 5.00$ we only find complex solutions to the system~\eqref{eq038} and~\eqref{eq040}. 
\begin{table}[ht]%
\begin{center}
\begin{tabular}{rrlrrr}
        &         &         & \multicolumn{3}{c}{\textbf{INITIAL VALUES}} \\
\cmidrule{4-6}        & \multicolumn{1}{l}{Value of the fixed parameter} & Curve   & $\ze_0$ & $\nu_{20}$ & $x_0$ \\
\midrule
\multicolumn{1}{l}{$\nu_1$: fixed} & \multicolumn{1}{l}{$\nu_1 =0.00$} & C1      & 1.00000 & 0.00000 & -1.00000 \\
        &         & C2      & 1.00000 & 0.19160 & -2.76929 \\
        &         & C3      & 4.00000 & 1.80565 & -2.22076 \\
\cmidrule{2-6}        & \multicolumn{1}{l}{$\nu_1 =0.05$} & C1      & 0.95000 & 0.04736 & -1.05573 \\
        &         & C2      & 0.95000 & 0.19197 & -2.72450 \\
        &         & C3      & 5.20000 & 0.79903 & -4.84245 \\
        &         & C4      & 5.20000 & 1.93559 & -9.57166 \\
\cmidrule{2-6}        & \multicolumn{1}{l}{$\nu_1 =1.50$} & C1      & 4.00000 & 0.03207 & -2.61688 \\
        &         & C2      & 4.00000 & 0.06738 & -6.74685 \\
\midrule
        &         &         & \multicolumn{3}{c}{} \\
        & \multicolumn{1}{l}{Value of the fixed parameter} &         & $\nu_{10}$ & $\nu_{20}$ & $x_0$ \\
\midrule
\multicolumn{1}{l}{$\ze$: fixed} & \multicolumn{1}{l}{$\ze =0.90$} & C1      & 0.00000 & 0.17053 & -1.14479 \\
        &         & C2      & 0.00000 & 0.22155 & -2.26820 \\
\cmidrule{2-6}        & \multicolumn{1}{l}{$\ze =5.00$} & C1      & 0.00000 & 0.03443 & -16.97040 \\
        &         & C2      & 0.00000 & 1.32591 & -2.83001 \\
        &         & C3      & 0.05200 & 2.06935 & -9.20408 \\
        &         & C4      & 1.20000 & 0.01599 & -3.10458 \\
        &         & C5      & 1.20000 & 0.05344 & -8.61426 \\
\cmidrule{2-6}        & \multicolumn{1}{l}{$\ze =8.00$} & C1      & 0.00000 & 0.02148 & -27.25100 \\
        &         & C2      & 0.00000 & 0.76381 & -4.62400 \\
        &         & C3      & 0.03400 & 0.56548 & -13.0963 \\
        &         & C4      & 0.70000 & 0.00952 & -5.03922 \\
        &         & C5      & 0.70000 & 0.03329 & -13.89100 \\
\bottomrule
\end{tabular}%

\caption{Initial conditions used to the computation of the critical damping curves shown in Fig.~\ref{fig02}}	
\end{center}
\label{tab01}
\end{table}
\begin{figure}[htb]
\begin{center}			
			\begin{tabular}{ccc}		
								Overdamped regions for $\ze \equiv$const & \phantom{aa} & Overdamped regions for $\nu_1 \equiv$const \\\\ \\
								\input{figures/figure02-S4.tex}  & \phantom{aa} &  	\input{figures/figure02-S1.tex} \\
								\input{figures/figure02-S5.tex}  &  					  &   \input{figures/figure02-S2.tex} \\
								\input{figures/figure02-S6.tex}  &  					  &  	\input{figures/figure02-S3.tex} \\ \\ \\
			\end{tabular}
			\caption{Example 2: Single degree--of freedom nonviscous system with $N=2$ kernels. Critical damping curves and overdamped regions}
			\label{fig02}
\end{center}
\end{figure}	
\newpage

In Fig.~\ref{fig02}, the different obtained critical curves have been plotted. If we read overdamped regions as certain volumes enclosed by critical surfaces in the parametric space $(\ze,\nu_1,\nu_2)$, then the left curves are cross sections of these volumes for certain values of the damping ratios (in Fig.~\ref{fig02} left cases $\ze=0.9,\ 5.0,\ 8.0$ are shown). On the other hand, right plots have the same interpretation but as cross sections of the planes $\nu_1=0.0,\ 0.054, \ 1.50$. The critical curves are in correspondence with the notation used in Table~\ref{tab01}. As expected overdamped region in the plane $(\nu_1,\nu_2)$ are symmetric respect to line $\nu_1=\nu_2$ since the physical model has inherently this symmetry. In the plane $(\nu_2,\ze)$ for $\nu_1=0$ we observe an overdamped region for values $\nu_2 \ll 1$ very similar to that obtained in the first example, see Fig.~\ref{fig01}. Now, due to the commented symmetry of the problem respect to $\nu_1$ and $\nu_2$, this long--triangle--like of the top--right plot is in fact a section of a conoid--like form in the space $(\ze,\nu_1,\nu_2)$. Actually, it is a quarter of conoid because $\nu_1,\nu_2 \geq 0$. Another section of this conoid is obtained just shifting the cross section to the value $\nu_1=0.05$ (middle--right plot). It is also interesting the new overdamped region arising in the corner of the right--top plot (section $\nu_1=0$). Presumably high values of $\ze$ and $\nu_2$ (simultaneously) lead the system to non--oscillatory motion. Let us see that the damper model in this case is formed by two dampers in parallel, one viscous and the other one nonviscous with parameter $\mu_1$. Indeed, for $\nu_1=0$ (which is equivalent to $\mu_1 \to \infty$) the damping function is transformed into
\begin{equation}
\lim_{\mu_1 \to \infty} \mG(t) =   c_1 \, \dt(t) +  c_2 \, \mu_2 \,  e^{-\mu_2 t}  \quad
\lim_{\mu_1 \to \infty} G(s) = c_1  + \frac{\mu_2 \, c_2}{s + \mu_2} 
\label{eq043}
\end{equation}
Somehow, we can interpret this overdamped region as the effect produced in the nature of the response of both the nonviscous parameter $\mu_2$ and the viscous coefficcient $c_1$. This is the reason because such critical surface did not appear in the example 1. Again, from the symmetry respect to the nonviscous parameters, this overdamped region also appears in the the plane $(\ze,\nu_1)$ for $\nu_2=0$. Furthermore, this effect is extended as an narrow volume in the approximate range $0 \leq \nu_1 \leq 0.07$ (respectively for symmetry in $0 \leq \nu_2 \leq 0.05$), see middle--left and bottom--left plots. It seems clear that adding new damping parameters will make more difficult how to read into the form of the overdamped manifolds, specially since they do not follow regular geometrical structures, as seen in this example. However, the proposed method could be applied sequentially to extract those more interesting curves for our analysis, for instance in artificial dampers design problems. Let us see now in a final example how to extract critical damping curves for a multiple dof system.

\subsection{Multiple degree of freedom systems}

\begin{figure}[H]%
\begin{center}
	\includegraphics[width=13cm]{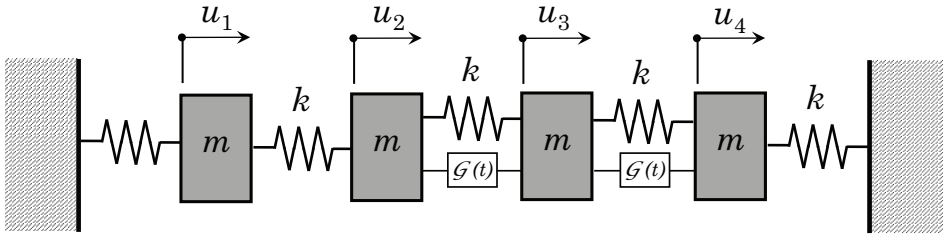} \\	
\caption{Example 3: The four degrees--of--freedom discrete system. $\mG(t)$ represents the hereditary function of nonviscous dampers}%
\label{fig04}%
\end{center}
\end{figure}

\begin{figure}[H]%
\begin{center}
	\input{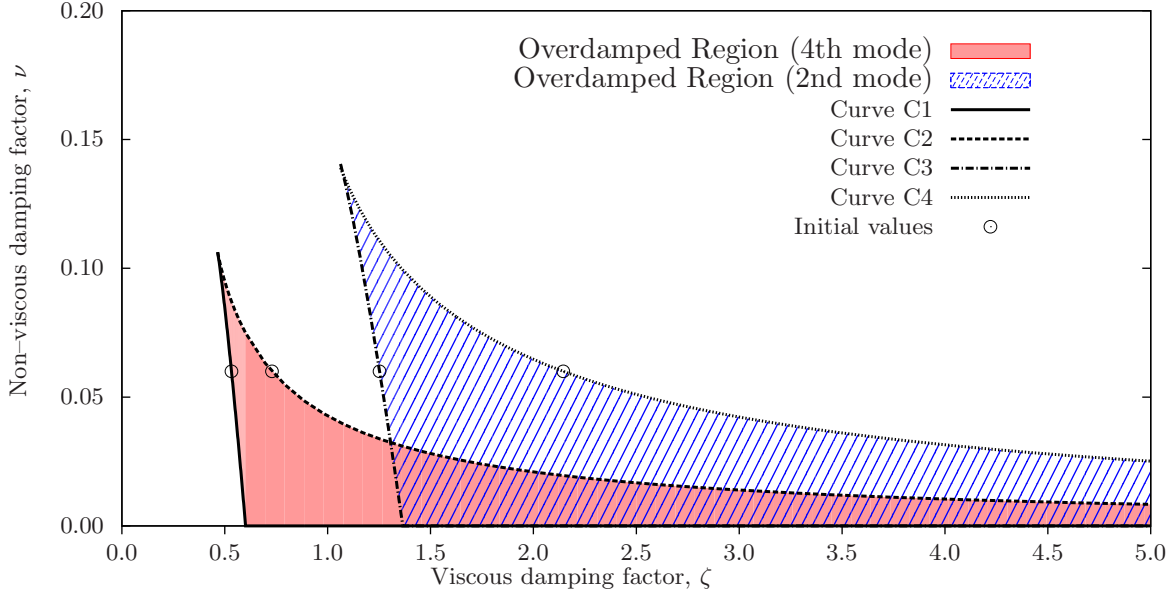}  \\	
\caption{Example 3: Critical damping curves and the  corresponding enclosing overdamped regions for the four degrees--of--freedom system}%
\label{fig03}%
\end{center}
\end{figure}

In order to validate the proposed approach to find critical damping curves for multiple dof systems a discrete lumped  mass dynamical system with four dof is analyzed. The Fig.~\ref{fig04} represents the distribution of masses $m$, rigidities $k$ and viscoelastic dampers with a hereditary  function $\mathcal{G}(t)$. The mass matrix of the system is $\mbf{M} = m \mbf{I}_4$ while, according to the rigidities and dampers distribution, stiffness matrix yields
\begin{equation}
		\mbf{K} = k \,
		\cor{\begin{array}{rrrr}
					2 	& -1 & 0 		& 0 \\
				- 1		& 2 & -1 & 0 \\
				   0 		& -1 & 2 & -1 \\
					 0 		&		 0 &	-1 & 2 \\
		\end{array}} = k \, \bmK
		\label{eq044}
\end{equation}								
We will assume a damping function formed by one hereditary exponential kernel of damping coefficient $c$ and nonviscous parameter $\mu$. Hence, te damping matrix can be expressed as 
$	\bmG(t) = \mu \, \mbf{C} \, e^{-\mu\,t}$ where
\begin{equation}
		\mbf{C} = c \,
		\cor{\begin{array}{rrrr}
					0 	& 	0 	& 	0 	& 0 \\
				  0		&		1		&	 -1		&	0 \\
				  0		&	 -1		&	  2		&	-1 \\
				  0		&	  0		&	  1		&	 1 \\
		\end{array}} \equiv c \, \bmC
		\label{eq045}
\end{equation}								
With help of these dimensionless matrices, say $\bmK$ and $\bmC$, we can express the non--linear eigenvalue problem associated to this problem under dimensionless form as
\begin{equation}
\cor{ x^2\mbf{I}_4 + \frac{2 \, x \, \ze}{1 + \nu \, x} \bmC + \bmK} \mbu = \mbf{0}
\label{eq046}
\end{equation}
where
\begin{equation}
x = \frac{s}{\om_0} \ , \quad \nu = \frac{\om_0}{\mu} \ , \quad \ze = \frac{c}{2 m \om_0} \ , \quad  \om_0 = \sqrt{k/m}
\label{eq047}
\end{equation}
Since $r=\text{rank} (\bmC) = 2$, the system has $2n+r=10$ eigenvalues. Therefore, the determinant of the transcendental matrix can be transformed into a 10th order polynomial multiplying by the factor $(1 + \nu x)^2$. We define then our function $\mD(x,\nu,\ze)$ as
\begin{eqnarray}
\mD(x,\nu,\ze) &=& (1 + \nu x)^2 \, \det \cor{x^2\mbf{I}_4 + \frac{2 \, x \, \ze}{1 + \nu \, x} \bmC + \bmK}  \nonumber  \\
							&=& 5 + 2 (8 \zeta +5 \nu ) x + \cor{ (2 \zeta +\nu ) (6 \zeta +5 \nu )+20}x^2 + \p{54 \zeta +40 \nu}x^3  \nonumber  \\
							& & + \p{32 \zeta ^2+54 \zeta  \nu +20 \nu ^2+21}x^4 + \p{40 \zeta +42 \nu}x^5 + \p{12 \zeta ^2+40 \zeta  \nu +21 \nu ^2+8} x^6  \nonumber  \\
							& &  + 8 (\zeta +2 \nu )x^7 + \p{1 + 8 \nu  \zeta +8 \nu^2 }x^8 + 2 \nu x^9 + \nu x^{10}
\label{eq048}
\end{eqnarray}
We are looking for overdamped regions enclosed by critical curves of type $\nu = \nu(\ze)$, therefore we construct our system of differential equations following the methodology described in Section~\ref{ProposedMethod}, resulting
\begin{equation}
\mmatriz{0}{\mD_{,\nu}}{\mD_{,xx}}{\mD_{,x\nu}} 
\ccolumna{x'(\ze)}{\nu'(\ze)} = 
- \ccolumna{\mD_{,\ze}}{\mD_{,x \ze}}
\ , \quad
\ccolumna{x(\ze_{0})}{\nu(\ze_{0})} = \ccolumna{x_0}{\nu_{0} } 
\label{eq049}
\end{equation}
For a sake of clarity in the exposition the expressions of the partial derivatives will not be written. Initial conditions can be found solving the system of algebraical equations for a particular value of $\nu_0$ and $\ze_0$
\begin{equation}
\mD(x_0,\nu_0,\ze_0)=0 \quad , \quad \mD_{,x}(x_0,\nu_0,\ze_0)=0
\label{eq050}
\end{equation}
Testing for $\nu_0=0.06$ we obtain four different pairs $(\ze_0,x_0)$, listed in Table~\ref{tab02}. 
\begin{table}[ht]%
\begin{center}
\begin{tabular}{lrrr}
        & \multicolumn{3}{c}{\textbf{INITIAL VALUES}} \\
\cmidrule{2-4}Curve   & $\ze_0$ & $\nu_{0}$ & $x_0$ \\
\midrule
C1      & 0.53258 & 0.06000 & -2.05512 \\
C2      & 0.72949 & 0.06000 & -7.88861 \\
C3      & 1.25218 & 0.06000 & -1.45645 \\
C4      & 2.14493 & 0.06000 & -8.07994 \\
\bottomrule
\end{tabular}%

\end{center}
\caption{Initial conditions used for the critical damping curves shown in Fig.~\ref{fig03}}	
\label{tab02}
\end{table}
After solving Eqs.~\eqref{eq049} we plot the solutions in Fig.~\ref{fig03}. The four found curves enclose two overdamped regions which in turn intersect each other. The solid--filled region represents the set of values $\ze,\nu$ which lead the fourth mode to overdamping. On the other hand, lines--filled area corresponds to the overdamped region of the second mode. This can be checked following a root--locus plot varying parameters $\ze,\nu$ from undamping to overdamping. Moreover, 2nd and 4th mode are precisely those modes for which degrees of freedom linked to the viscoelastic dampers are most activated. Obviously, it follows then that the overlapping zone (with both types of shading in Fig.~\ref{fig03}) corresponds with the overdamping of both modes, simultaneously. 
A deeper inspection of Fig.~\ref{fig03} leads us to ask ourselves about the existence of singularities in the domain $(\ze,\nu)$ which can impede the application of our approach. Well, it is known that a system of differential equations like that shown in \eqref{eq049} has solution (and it is unique) provided that the determinant of the matrix does not vanish in a neighborhood around the initial value. This determinant is equal to $-\mD_{,\nu}\,\mD_{,xx}$, whence it follows that these two system of algebraic equations 
\begin{equation}
\text{S1:}
\begin{cases} \mD=0 \\ \mD_{,x}=0 \\ \mD_{,\nu} =0 \end{cases}
\quad \ \quad
\text{S2:}
\begin{cases} \mD=0 \\ \mD_{,x}=0 \\ \mD_{,xx}=0 \end{cases}
\label{eq051}
\end{equation}
allow us to find singularities. For every solution of the first system S1, some of the variables, $x$, $\ze$ or $\nu$ is within the complex plain. On the other hand, we do find valid solutions of system S2, verifying $x<0, \ \ze,\nu\geq 0$, say
\begin{align}
x& = -3.1059		& \ze &= 0.46473		& \nu &= 0.10658  \nonumber \\
x& = -2.3221		& \ze &= 1.06113		& \nu &= 0.14088 \label{eq052}  
\end{align}
These two points are precisely the vertexes of the two overdamped regions, namely, intersection points of curves C1--C2 and C3--C4, respectively. These points can not be used as initial points since according to the implicit function theorem, equations $\mD=0$ and $\mD_{,x}=0$ do not define $x(\ze)$ and $\nu(\ze)$ unequivocally. Furthermore, both points verifies $\mD_{,xx}=0$, hence they are triple roots. \\

We wonder now how does the solution behave in the intersection point between curves C2 and C3, located approximately at $\ze=1.30465, \ \nu=0.03239$. This point satisfies $-\mD_{,\nu}\,\mD_{,xx}\neq 0$, therefore it should be valid as initial value of our method. However, it belongs to both curves simultaneously, so that at a first sight the solution would not seem to be well defined. However, in this point we find two solutions for the variable $x$, say $x = -1.38042$ and $x=-15.2148$, which leads to two different initial values. Somehow this point does not correspond with a intersection point of the curves in the 3D domain $(x,\nu,\ze)$, while the {\em triple} roots of \eqref{eq052} do.\\

As far as multiple degree-of-freedom systems concern, the success of the method lies on the availability of the transcendental matrix determinant and their derivatives. Therefore, from a numerical point  of view, large systems will require high computational effort which limits the range of applicability for small or moderate order systems. Currently, our efforts are addressed to find out numerical procedures allowing to construct approximate critical curves but for larger systems, something that is under research.

\section{Conclusions}

In this paper critical damping of nonviscously damped linear systems is studied. Nonviscous or viscoelastic vibrating structures are characterized by dissipative mechanisms depending on the history of response through hereditary functions. For certain values of the damping parameters, the response can become non--oscillatory. It is said then that some (or all) modes  are overdamped. Particular values of the damping parameters which establish the limit between oscillatory and non--oscillatory motion are said to be within a critical surface (or critical manifold). In the present paper a general procedure to build critical damping surfaces is developed. In addition, a numerical method based on the transformation of the algebraical equations into a system of two ordinary differential equations is proposed. This approach allow to find critical curves of two parameters for certain fixed values of the rest of parameters. \\

To validate the theoretical results three numerical examples are analyzed. In the first example, the well--known overdamped region of a single degree--of--freedom system with one exponential kernel is resolved, shown perfect fitting between our curves, obtained from the differential equations, and those of the analytical expressions. In order to give  added value to this problem, we propose simplified approximate expressions for the critical curves. The second example is devoted to construct overdamped regions for a two exponential kernels based damping function. This problem involves three parameters, so that the critical damping curves are plotted along several cross sections defined by the free parameter. The third example shows how the method can be applied for multiple degrees of freedoms systems deriving overdamping regions for different modes and interpreting the obtained overlapping regions. Since this method is based on the evaluation and derivation of the determinant of the transcendental matrix, its range of validity is reduced to small or moderately sized systems. Encouraged by this limitation, the author is currently investigating how to extrapolate this method for larger systems.








\end{document}